\documentclass[twoside,11pt]{article}
\title{} \author{} \date{}
\usepackage{latexsym,amssymb,times}
\input amssym.def
\newtheorem{te}{Theorem}[section]
\newtheorem{fac}[te]{Fact}
\newtheorem{lem}[te]{Lemma}

\newtheorem{ex}[te]{Example}

\newtheorem{df}[te]{Definition}

\def\dok{\noindent{\bf Proof. }}
\def\kdok{\hfill $\Box$ \par \vspace*{2mm} }

\def\a{\alpha}

\def\f{\varphi}
\def\p{\psi}

\def\o{\omega}

\def\r{\rho}
\def\s{\sigma}

\def\T{{\mathbb T}}
\def\P{{\mathbb P}}
\def\Q{{\mathbb Q}}
\def\B{{\mathbb B}}
\def\N{{\mathbb N}}

\def\X{{\mathbb X}}
\def\Y{{\mathbb Y}}
\def\Z{{\mathbb Z}}
\def\D{{\mathbb D}}
\def\C{{\mathbb C}}
\def\A{{\mathbb A}}
\def\BG{{\mathbb G}}

\def\R{{\mathbb R}}

\def\CP{{\mathcal P}}

\def\CA{{\mathcal A}}

\def\la{\langle}
\def\ra{\rangle}

\def\up{\!\upharpoonright\!}

\def\akko{\Leftrightarrow}
\def\dom{\mathop{\mathrm{dom}}\nolimits}
\def\ran{\mathop{\mathrm{ran}}\nolimits}
\def\id{\mathop{\mbox{id}}\nolimits}
\def\Emb{\mathop{\rm Emb}\nolimits}

\def\Aut{\mathop{\rm Aut}\nolimits}

\def\Pi{\mathop{\rm Pi}\nolimits}

\def\ar{\mathop{\mathrm{ar}}\nolimits}

\def\orb{\mathop{\mathrm{orb}}\nolimits}
\def\Age{\mathop{\mathrm{Age}}\nolimits}

\def\supp{\mathop{\mathrm{supp}}\nolimits}
\def\ps{\mathbin{\emptyset}}
\newcommand{\seq}[1]{\left<#1\right>}
\newcommand{\set}[1]{\left\{#1\right\}}
\newcommand{\abs}[1]{\left\vert#1\right\vert}
\begin{document}
\thispagestyle{plain}
\begin{center}
          {\large \bf
                    \uppercase{Antichains of copies of ultrahomogeneous structures}}
\end{center}
\begin{center}
{\bf Milo\v s S.\ Kurili\'c\footnote{Department of Mathematics and Informatics, University of Novi Sad,
     Trg Dositeja Obradovi\'ca 4, 21000 Novi Sad, Serbia. e-mail: milos@dmi.uns.ac.rs}
     and
     Bori\v sa  Kuzeljevi\'c\footnote{Department of Mathematics and Informatics, University of Novi Sad,
     Trg Dositeja Obradovi\'ca 4, 21000 Novi Sad, Serbia. e-mail: borisha@dmi.uns.ac.rs}
}
\end{center}
\begin{abstract}
\noindent
We investigate possible cardinalities of maximal antichains in the poset of copies $\la \mathbb P(\mathbb X),\subset \ra$
of a countable ultrahomogeneous relational structure $\X$.
It turns out that if the age of $\X$ has the strong amalgamation property, then,
defining a copy of $\X$ to be large iff it has infinite intersection with each orbit of $\X$,
the structure $\X$ can be partitioned into countably many large copies,
there are almost disjoint families of large copies of size continuum
and, hence, there are (maximal) antichains of size continuum in the poset $\P (\X)$.
Finally, we show that the posets of copies of all countable ultrahomogeneous partial orders contain maximal antichains of cardinality continuum
and determine which of them contain countable maximal antichains. That holds, in particular, for the random (universal ultrahomogeneous) poset.


\end{abstract}
\section{Introduction}\label{S1}
In this paper we investigate antichains in the posets of the form $\la \mathbb P(\mathbb X),\subset \ra$,
where $\P (\X ):= \{ f[X]\!:\!f\!\in \Emb (\X)\}$ is the set of the substructures of a countable ultrahomogeneous relational structure $\X$
which are isomorphic to $\X$.
Recall that a structure $\X$ is \emph{ultrahomogeneous} iff for each isomorphism $\varphi:\mathbb A\to \mathbb B$
between finite substructures $\A$ and $\B$ of $\X$,
there is an automorphism $f$ of $\X$ extending $\varphi$.
These posets were analyzed from various viewpoints recently.
Typically, the results obtained would be compared to the poset $\la [\omega]^{\omega},\subset \ra$
of all infinite subsets of a countable set, ordered by inclusion.
Set theorists thoroughly investigated this object
and, most often, an antichain in this context is a set of pairwise incompatible elements,
i.e.\ a collection of sets in $[\o]^{\o}$ with pairwise finite intersections (an {\it almost disjoint family}).
Two basic facts are that there is no countable maximal antichain in $[\o]^{\o}$, whereas there is a maximal antichain of size continuum in that poset.
We follow this approach.
So, in this paper, an antichain is always a set of pairwise incompatible elements of the partial order in question.

Section \ref{S2} contains definitions and facts which are used in the paper. Defining a copy of $\X$ to be {\it large}
iff it has infinite intersection with each orbit of $\X$, in Sections \ref{dc} and \ref{mainsec} we prove the following general statement.
\begin{te}\label{T9050}
If $\,\X$ is a countable ultrahomogeneous relational structure satisfying the strong amalgamation property, then\\[-7mm]
\begin{itemize}\itemsep=-1.5mm
\item[(a)] $\X$ can be partitioned into countably many large copies of $\,\X$;
\item[(b)] There are almost disjoint families of large copies of $\,\X$ of size continuum;
\item[(c)] There are (maximal) antichains of size continuum in the poset $\la \mathbb P(\mathbb X),\subset \ra$.
\end{itemize}
\end{te}
In Section \ref{posets} we take a closer look on the case of countable ultrahomogeneous posets, using the following
well-known classification due to Schmerl \cite{Sch}.
\begin{te}[Schmerl]\label{schmerl}
Each countable ultrahomogeneous partial order is isomorphic to one of the following:\\[-7mm]
\begin{itemize}\itemsep=-1.5mm
\item[-] $\mathbb A _\omega$, a countable antichain (that is, the empty relation on $\omega$);
\item[-] $\mathbb B _n = n \times \mathbb Q$, $1\leq n\leq \omega$, where $\langle i_1, q_1\rangle < \langle i_2 ,q_2\rangle \Leftrightarrow i_1=i_2 \;\land \; q_1 <_\mathbb Q q_2$;
\item[-] $\mathbb C _n=n\times \mathbb Q$, $1\leq n\leq \omega$, where  $\langle i_1, q_1\rangle < \langle i_2 ,q_2\rangle \Leftrightarrow  q_1 <_\mathbb Q q_2$;
\item[-] $\mathbb D$, the random poset.
 \end{itemize}
\end{te}
So, in Section \ref{posets} we prove the following theorem.
%
\begin{te}\label{T9092}
Let $\X$ be a countable ultrahomogeneous partial order. Then\\[-7mm]
\begin{itemize}\itemsep=-1.5mm
 \item[(a)] There are maximal antichains of size continuum in the poset $\la \mathbb P(\mathbb X),\subset \ra$;
 \item[(b)] There are countable maximal antichains in $\la \mathbb P(\mathbb X),\subset \ra$ if and only if
            $\X$ is neither isomorphic to $\mathbb A_{\o}$ nor to $\mathbb B_{\omega}$.
 \end{itemize}
\end{te}
At this point we mention some related concepts.
First, antichains in the poset of copies of the random (Rado) graph were analyzed in \cite{KP}.
Second, forcing-related properties of the posets of copies of ultrahomogeneous structures were investigated in \cite{KTQ,KTR,KTR1}.
Third,  in \cite{Ktow,Kurscatt,Kurord}  a classification of relational structures
with respect to the properties of posets $\la \mathbb P(\mathbb X),\subset \ra$ is given.
Fourth,  the order types of the maximal chains in the posets of copies
of countable ultrahomogeneous graphs and countable ultrahomogeneous partial orders
are described in \cite{KK1,KK2}.
Finally, if $\X$ is a first order structure and $\preceq ^R$ right Green's pre-order on its self-embedding monoid, $\Emb \X $,
the corresponding antisymmetric quotient $\la \Emb \X /\!\!\approx ^R, \preceq ^R\ra$ (right Green's order)
is isomorphic to the partial order $\la \P (\X ), \supset \ra$.
Hence, our results provide some information about self-embedding monoids of structures.
\section{Preliminaries}\label{S2}
If $L=\la R_i : i\in I\ra$ is a relational language, where $\ar (R_i)=n_i \in \N$, for $i\in I$, and
$\X =\la X, \r \ra$ is an $L$-structure, where  $\r=\la \r _i :i\in I\ra$ and $\r _i \subset X^{n_i}$, for $i\in I$, then, for a subset $A$ of $X$, by
$\r \up A$ we will denote the sequence $\la \r _i \up A :i\in I\ra$, where  $\r _i \up A := \r _i \cap A^{n_i}$. Then the structure $\A=\la A , \r \up A\ra$ is a substructure of $\X$.

If $\Y =\la Y, \s \ra$ is also an $L$-structure, an injection
$f:X \rightarrow Y$ is called an {\it embedding} (we write $f: \X \hookrightarrow  \Y$ or $f\in \Emb (\X , \Y )$) iff for each $i\in I$
and $\bar x \in X^{n_i}$ we have: $\bar x \in \r _i $ iff $f\bar x \in \s _i$.
If, in addition, $f$ is a surjection, it is an {\it  isomorphism}, the structures $\X$ and $\Y$ are {\it  isomorphic},
and we write $\X \cong \Y$.
If, in particular, $ \Y = \X $, then $f$ is an {\it  automorphism} of the structure $\X$.
$\Aut (\X )$ will denote the set of all automorphisms of $\X$ and
$\Emb (\X )$ is $\Emb (\X , \X ) $.
By $\P (\X )$ we denote the set of domains of substructures of $\X$ isomorphic to $\X$, that is
$\P (\X )  =  \{ A\subset X : \la A, \r \up A \ra \cong \la X, \r \ra\}=  \{ f[X] : f \in \Emb (\X )\} $.

If $\seq{P,\leq }$ is a poset, the elements $x$ and $y$ of $P$ are \emph{compatible} iff there is an element $z\in P$ such that $z\le x$ and $z\le y$.
Otherwise, $x$ and $y$ are {\it incompatible} and we write $x\perp y$.
A set $A\subset P$ is \emph{an antichain} in $P$ if its elements are pairwise incompatible.
An antichain $A$ is a \emph{maximal antichain} in $P$ iff each $z\in P$ is compatible with some $x\in A$.

We recall some basic facts from Fra\"{\i}ss\'{e} theory.
The {\it age}, $\Age \X$, of an ultrahomogeneous $L$-structure $\X$ (i.e., the class  of all
finite $L$-structures embeddable in $\X$)
satisfies
the {\it amalgamation property} (AP): if $\A ,\B ,\C \in \Age \X$ and $f_0:\A\hookrightarrow \B$ and $g_0:\A\hookrightarrow \C$ are embeddings,
then there are $\D\in \Age \X$ and embeddings
$f_1:\B\hookrightarrow \D$ and $g_1:\C\hookrightarrow \D$ such that $f_1 \circ f_0=g_1 \circ g_0$.
If, in addition, the amalgam $\D$ and the embeddings $f_1$ and $g_1$ can be chosen so that $f_1[B]\cap g_1[C]=f_1[f_0[A]]=g_1[g_0[A]]$,
then (the age of) $\X$ satisfies the {\it strong amalgamation property} (SAP).
We will use the following classical results of Fra\"{\i}ss\'e (see \cite{Fra}, p.\ 332--333).
\begin{te}\label{T9006}
(a) A countable structure $\X$ is ultrahomogeneous iff for each finite substructure $\A$ of $\X$, each $f\in\Emb(\A,\X)$ and each $x\in X\setminus A$
there is $y\in X$ such that $f\cup\set{\seq{x,y}}\in \Emb(\A\cup\set{x},\X)$.

(b) Countable ultrahomogeneous structures with the same age are isomorphic.
\end{te}
If $\X$ is an $L$-structure, the {\it pointwise stabilizer} of  a finite set $F\subset X$ is the subgroup
$\Aut _F (\X ):=\{ g\in \Aut (\X ): \forall x\in F\; g(x)=x \}$ of the group $\Aut (\X )$.
The binary relation $\sim _F$ on the set $X\setminus F$ defined by
$x\sim _F y $ iff there is $g\in \Aut _F (\X )$ such that $g(x)=y$,
is an equivalence relation and the equivalence class of an $x\in X\setminus F$ is
denoted by $\orb _{F}(x)$ and called  the {\it orbit of $x$ under $\Aut _F(\X )$}. Thus
$$
\orb _F(x)=\Big\{ y\in X\setminus F :\exists g\in \Aut (\X )\;\; ( g\up F =\id _F \land \;g(x)=y )\Big\}.
$$
The sets $\orb _F (x)$, where $F\in [X ]^{<\o}$ and $x\in X\setminus F$, are called {\it the orbits of} $\X$.
Later in the paper, the strong amalgamation property will play a significant role and the next theorem provides convenient characterizations of this property.
\begin{te}\label{T9044}(see \cite{Fra} p.\ 399 and \cite{Camer} p.\ 37)
For a countable ultrahomogeneous relational structure $\X$  the following conditions are equivalent:\\[-7mm]
\begin{itemize}\itemsep=-1.5mm
\item[(a)] $\X$ satisfies the strong amalgamation property,
\item[(b)] $\X$ is strongly inexhaustible, that is, $X\setminus F \in \P (\X)$, for each finite $F\subset X$,
\item[(c)] The orbits of $\X$ are infinite.
\end{itemize}
\end{te}
The following characterization of copies of ultrahomogeneous structures is, most likely, a known fact. We include its proof for completeness of the paper.
\begin{te}\label{T9086}
If $\X$ is a countable ultrahomogeneous $L$-structure and $A\subset X$, then $A\in \P(\X)$ iff
\begin{equation}\label{EQ9099}
\forall F\in [A]^{<\o}\ \forall x\in X\setminus F\ \orb_F(x)\cap A\neq\ps.
\end{equation}
Consequently, if the set $A$ intersects all orbits of $\X$, then $A\in \P(\X)$.
\end{te}
\dok
Let $\A\cong\X$, $F\in [A]^{<\o}$ and $x\in \X\setminus F$.
Then there are $F'\subset A$, $a'\in A$ and an isomorphism $\psi:F'\cup\{a'\} \rightarrow F\cup \{x\}$ such that $\psi (a')=x$.
Since the structure $\A$ is ultrahomogeneous, by Theorem \ref{T9006}(a) there is $a\in A$ such that
$\f :=(\psi\up F') \cup \{ \la a',a\ra\}: F'\cup\{a'\} \rightarrow F\cup \{a\}$ is an isomorphism.
Now $\eta:=\varphi\circ \psi ^{-1}$ is an isomorphism, $\eta\up F =\id _F$ and $\eta (x)=a$. Since $\X$ is ultrahomogeneous,
there is $g\in \Aut (\X )$ extending $\eta$; so $g\in \Aut _F(\X )$, $g(x)=a$ and $a\in \orb_F(x)\cap A$.

Assuming (\ref{EQ9099}) we prove that the set $\Pi (\A,\X) $ of all finite partial isomorphisms from $\A$ into $\X$
has the back-and-forth property.  So, let $\f \in \Pi (\A,\X)$.
First, if $a\in A\setminus \dom \f$, then
by Theorem \ref{T9006}(a) there is $x\in X$ such that $\p := \f \cup \{\la a ,x\ra\}$ is an isomorphism and, clearly,
$\p \in \Pi (\A,\X)$.
Second, if $x\in X \setminus \ran \f$, then
by Theorem \ref{T9006}(a) there is $x'\in X$ such that $\p := \f ^{-1}\cup \{\la x ,x'\ra\}$ is an isomorphism.
Since $x'\not\in \dom \f$, by (\ref{EQ9099}) there exists $a\in \orb _{\dom \f}(x')\cap A$
and, hence, there is $g\in \Aut _{\dom \f}(\X)$ such that $g(x')=a$.
Now, $\eta := g\circ \p$ is a finite isomorphism with domain $\ran \f \cup \{ x \}$ and
$\eta [\ran \f \cup \{ x \}]=g[\p [\ran \f]]\cup g[\p [\{ x \}]]=\dom \f \cup \{ a \}$. In addition, for $y\in \ran \f$
we have $\eta (y)=g(\p (y))=g(\f ^{-1}(y))=\f ^{-1}(y)$, which gives $\eta =\f ^{-1}\cup \{ \la x,a \ra\}$.
So $\eta ^{-1} =\f \cup \{ \la a,x \ra\}\in \Pi (\A,\X)$. 
\hfill $\Box$
\section{Partitions into large copies}\label{dc}
Here we make some observations about copies of ultrahomogeneous structures incompatible in a very strong way.
By $DC$ we denote the class of countable structures {\it having disjoint copies}
(there are  copies $A,B\in \P (\X )$ such that $A\cap B=\emptyset$)
and by $SAP$ the class of countable ultrahomogeneous structures satisfying SAP.
\begin{ex}\label{EX9002}\rm
{\it A countable ultrahomogeneous structure without disjoint copies.}
Let $\X =\Q \cup \Y$, where $\Q$ is the rational line and $\Y :=\la \{ y \} , \{ \la y,y \ra\}\ra$,
where $y\not \in \Q$. Then $y\in A$, for each $A\in \P (\X )$.
\end{ex}
A structure $\X $ is called {\it indivisible} (resp.\ {\it strongly indivisible}) iff for each partition
$X=A \cup B$ there is $C\in \P (\X )$ such that $C\subset A$ or $C\subset B$ (resp.\ $A\in \P (\X )$ or $B\in \P (\X )$).
Let $UH$, $I$, and $SI$, denote the classes of ultrahomogeneous, indivisible and strongly indivisible  countable relational structures respectively.

Confirming a conjecture of Fra\"{\i}ss\'{e},
Pouzet proved that each countable indivisible structure $\X$ has disjoint copies \cite{Pouz}; thus, $I \subset DC$.
Here we prove that more holds for countable ultrahomogeneous structures satisfying SAP; thus $SAP \subset DC$.
We note that $SAP \not\subset I$, for example, $\B _n \in SAP \setminus I$, for $1<n<\o$.
\begin{te}\label{T9048}
Each countable ultrahomogeneous structure $\X$ satisfying SAP can be partitioned into countably many large copies of $\X$.
\end{te}
\dok
W.l.o.g.\ we assume that $X=\o$. By Theorem \ref{T9044}, the set of orbits, 
$\Omega := \{ \orb _F (x): F\in [\o ]^{<\o } \land x\in \o\setminus F\}$,
is a countable subfamily of $[\o ]^{\o }$.
Let $\Omega =\{ O_n : n\in \o \}$ be an enumeration of $\Omega$
and let the sequence $\la m_{n,k} :n\leq k< \o \ra$ in $\o $ be constructed by recursion as follows. First, let $m_{0,0}=\min O_0$.

Second, if $0<k<\o$ and $m_{n',k'}$ are defined for $\,n'\leq k' <k$, then we define $m_{n,k}$ for $n\leq k$ by:
$m_{0,k}=\min[ O_0\setminus \{ m_{n',k'} : n'\leq k' <k\}]$ and, for $0< n\leq k$,
\begin{equation}\label{EQ9006}
m_{n,k}=\min\Big[ O_n\setminus \Big(\{ m_{n',k'} : n'\leq k' <k\} \cup \{ m_{n',k} : n' <n\}\Big)\Big].
\end{equation}
\noindent
The recursion works, since $|O_n|=\o$, for all $n\in \o$.
By the construction, all the $m_{n,k}$'s are different. So, defining $A_i:=\{ m_{n,n+i}:n<\o \}$, for all $i\in \o$,
we have $A_{i_1}\cap A_{i_2} =\emptyset$, for $i_1\neq i_2$.
By (\ref{EQ9006}), for each $n\in \o$ we have $m_{n,n+i}\in A_i\cap O_n$,
thus the set $A_i$ intersects all the orbits of $\X$ and, by Theorem \ref{T9086}, $A_i\in \P (\X )$.

Clearly we have $A:=\bigcup _{i\in \o}A_i =\{ m_{n,k} :n\leq k< \o \}\subset \bigcup _{n\in \o}O_n=X$.
Suppose that $O_n\not\subset A$, for some $n\in \o$ and let $m=\min (O_{n} \setminus A)$.
By (\ref{EQ9006}) we have $|A\cap O_n|=|\{ m_{n,k} : k< \o\}|=\o$, so there is $k=\min \{ k'<\o :m<m_{n,k'}\}$
and we have $m< m_{n,k}$.
Since $\{ m_{n',k'} : n'\leq k' <k \} \cup \{ m_{n',k} : n' <n\}\subset A$
we have $S:= O_{n} \setminus\{ m_{n',k'} : n'\leq k' <k \} \cup \{ m_{n',k} : n' <n\}\supset O_{n} \setminus A $
and, by  (\ref{EQ9006}), $m_{n,k}=\min S\leq \min (O_{n} \setminus A)=m$, which gives a contradiction.
Thus $A=X$ and $\{ A_i :i<\o\}$ is a partition of $X(=\o )$.

\begin{figure}[h]
\begin{center}
\unitlength 0.8mm 
\linethickness{0.5pt}
\ifx\plotpoint\undefined\newsavebox{\plotpoint}\fi 
\begin{picture}(110,90)(0,0)
\put(5,5){\line(1,0){100}}
\put(25,10){\line(1,0){40}}
\put(10,15){\line(1,0){90}}
\put(30,20){\line(1,0){55}}
\put(30,40){\line(1,0){55}}
\put(10,50){\line(1,0){70}}
\put(90,50){\line(1,0){10}}
\put(25,60){\line(1,0){5}}
\put(40,60){\line(1,0){25}}
\put(5,65){\line(1,0){100}}
\put(5,85){\line(1,0){25}}
\put(40,85){\line(1,0){25}}
\put(5,5){\line(0,1){80}}
\put(10,15){\line(0,1){35}}
\put(25,10){\line(0,1){50}}
\put(30,20){\line(0,1){20}}
\put(65,5){\line(0,1){80}}
\put(85,20){\line(0,1){5}}
\put(85,35){\line(0,1){5}}
\put(100,15){\line(0,1){35}}
\put(105,5){\line(0,1){25}}
\put(105,40){\line(0,1){25}}



\footnotesize
\put(17,30){\makebox(0,0)[cc]{$\X _{Lach}$}}
\put(45,30){\makebox(0,0)[cc]{$\BG _{Rado}  $}}
\put(75,30){\makebox(0,0)[cc]{$\la \o , <\ra  $}}
\put(85,30){\makebox(0,0)[cc]{$SI  $}}
\put(105,35){\makebox(0,0)[cc]{$DC  $}}
\put(45,45){\makebox(0,0)[cc]{$\Q  $}}
\put(75,45){\makebox(0,0)[cc]{$[0,1]_\Q  $}}
\put(85,50){\makebox(0,0)[cc]{$I  $}}
\put(17,55){\makebox(0,0)[cc]{$\bigcup _\o T_3  $}}
\put(45,55){\makebox(0,0)[cc]{$\B _2  $}}
\put(75,55){\makebox(0,0)[cc]{$\la \Z ,< \ra  $}}
\put(35,60){\makebox(0,0)[cc]{$SAP  $}}
\put(17,75){\makebox(0,0)[cc]{$\Q \cup 1_{re}  $}}
\put(75,75){\makebox(0,0)[cc]{$ \BG _\o $}}
\put(35,85){\makebox(0,0)[cc]{$UH  $}}
\end{picture}
\end{center}
\vspace{-5mm}
\caption{Countable relational structures}\label{F4002}
\end{figure}
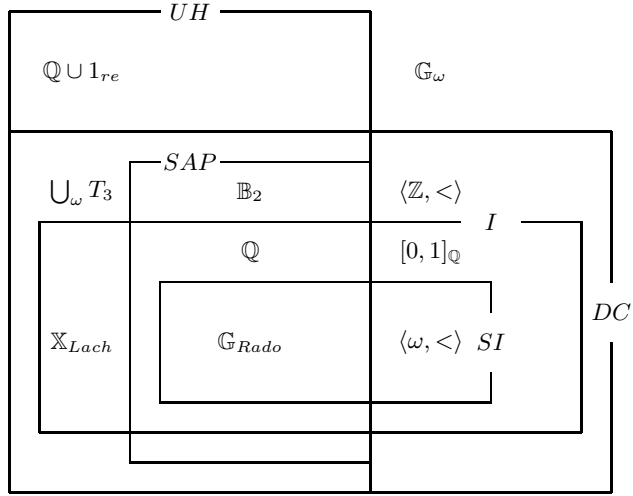

Now, let $\{ S_j :j\in \o\}\subset [\o ]^\o$ be a partition of $\o$ and $B_j:=\bigcup _{i\in S_j}A_i$, for $j\in \o$.
Then $\{ B_j :j\in \o\}\subset \P (\X )$ is a partition of $X$ and for $n\in \o$ we have
$\{ m_{n,n+i}:i\in S_j \}\subset B_j \cap O_n$; thus $B_j$, $j\in \o$, are large copies of $\X$.
\hfill $\Box$
\begin{ex}\label{EX9003}\rm
{\it A countable ultrahomogeneous divisible structure which does not have the SAP, but has disjoint copies.}
Let $\X$ be the wreath product $I_\o [T_3]$ (see \cite{Che1}), that is the disjoint union $\bigcup _{n\in \o} T_3^n$ of $\o$-many copies of the oriented triangle.
By Theorem \ref{T9044} the structure $\X$ does not satisfy SAP,
it is clear that $\X$ is not indivisible, but for each $S\in [\o ]^{\o}$ we have $A_S:=\bigcup _{n\in S} T_3^n \in \P (\X )$.
Concerning Theorems \ref{T9086} and \ref{T9048} we note that each one-element subset of $X$ is an orbit of $\X$.
Hence $X$ is the only subset of $X$ intersecting all the orbits of $\X$.
\end{ex}
Figure \ref{F4002} shows the relationship between the mentioned five classes. For $\X _{Lach}$ see \cite{Fra}, p.\ 402. $\mathbb G_{\omega}$ is the linear graph on $\omega$, i.e. $\seq{\o,\sim}$, where $m\sim n\akko \abs{m-n}=1$. $\mathbb Q\cup 1_{re}$ is the structure from Example \ref{EX9002}.
\section{Large almost disjoint families of large copies}\label{mainsec}
By Theorem 5.3 of \cite{Ktow}, if $\X$ is a countable indivisible $L$-structure, then the set $\P (\X )$ contains an almost disjoint family of size ${\mathfrak c}$
and, hence, the poset $\P (\X )$ contains maximal antichains of size ${\mathfrak c}$.
Here, proving Theorem \ref{T9050}(b) and (c), we show that more is true for countable ultrahomogeneous structures satisfying SAP.
\begin{te} \label{T9052'}
If $\,\X$ is a countable ultrahomogeneous $L$-structure satisfying SAP,
then there is an almost disjoint family of large copies of size ${\mathfrak c}$.
\end{te}
\dok
W.l.o.g.\ we assume that $X=\o$. By Theorem \ref{T9044}, the set of orbits, 
$\Omega := \{ \orb _F (x): F\in [\o ]^{<\o } \land x\in \o\setminus F\}$,
is a countable subfamily of $[\o ]^{\o }$.
Let $\Omega =\{ O_n : n\in \o \}$ be an enumeration of $\Omega$
and let  $\la m_{n,k} :n\leq k< \o \ra$  be the sequence in $\o $ constructed in Theorem \ref{T9048}.
Namely, $m_{0,0}=\min O_0$ and
if $m_{n',k'}$ are defined for $\,n'\leq k' <k$, then
\begin{equation}\label{EQ9006'}
m_{n,k}=\min\Big[ O_n\setminus \Big(\{ m_{n',k'} : n'\leq k' <k\} \cup \{ m_{n',k} : n' <n\}\Big)\Big].
\end{equation}
Since all the $m_{n,k}$'s are different, defining $D_n :=\{ m_{n,k} : k\in [n, \o)\} $, for $n\in \o$, we have
$D_n \in [O_n]^\o$ and $D_m\cap D_n =\emptyset$, for $m\neq n$.
Clearly we have $D:=\bigcup _{n\in \o}D_n =\{ m_{n,k} :n\leq k< \o \}$
and (see the proof of Theorem \ref{T9048}) $D=\o$. So $\{ D_n :n<\o\}$ is a partition of $\o$ refining $\Omega$.\\[-2mm]

\noindent
{\it Claim.}
If $\{ D_n :n<\o\}\subset [\o ]^\o$ is a partition of $\o$, then
there exists an almost disjoint family 
$\{ A_\a :\a <{\mathfrak c}\}\subset [\o ]^\o$, such that $|A_\a \cap D_n|=\o$, for each $\a <{\mathfrak c}$ and each $n\in \o$.\\[-2mm]

\noindent
{\it Proof of Claim.}
W.l.o.g.\  instead of $\o$ we can take the set of rationals, $\Q$, and suppose that $D_n$, $n<\o$, are dense suborders of $\Q$.
Let $f:\o \rightarrow \o$ be a surjection such that $|f^{-1}[\{ n \}]|=\o$, for each $n\in \o$.
By recursion, for each real $x\in \R$  we construct an increasing sequence $\la q^x_k:k\in \o\ra$ in $\Q$ converging to $x$
in the following way. First we take $q^x_0\in D_{f(0)}\cap (-\infty ,x)$; 
if $q^x_0 ,\dots ,q^x_k$ are defined, then we take $q^x_{k+1}\in D_{f(k+1)}\cap (\max \{q^x_k, x-\frac{1}{k+1}\} ,x)$. By the density of the sets $D_n$
the recursion works. Now, defining the sets $A_x :=\{ q^x_k:k\in \o \}$, for $x\in \R$, and ${\mathcal A}:=\{ A_x :x\in \R\}\in [\Q]^\o$
we obtain an almost disjoint family of size ${\mathfrak c}$. In addition, for $x\in \R$, $n\in \o$ and $k\in f^{-1}[\{ n \}]$
we have $q^x_{k}\in A_x \cap D_{f(k)}=A_x \cap D_n$ and, since $|f^{-1}[\{ n \}]|=\o$ and $q^x_{k}$'s are different, we have $|A_x \cap D_n|=\o$.
\kdok
By Claim, there is an almost disjoint family $\{ A_\a :\a <{\mathfrak c}\}\subset [\o ]^\o$
such that for each $\a <{\mathfrak c}$ and each $n\in \o$ we have $|A_\a \cap D_n|=\o$
and, since $D_n\subset O_n$, $|A_\a \cap O_n|=\o$. By Theorem \ref{T9086} we have $A_\a \in \P (\X )$.
\hfill $\Box$
\begin{ex}\rm \label{EX9000}
{\it Applications of Theorem \ref{T9050}.} The countable ultrahomogeneous digraphs
(structures with one irreflexive and asymmetric binary relation) have been classified by Cherlin \cite{Che1,Che2}.
Following the organization of the {\it Cherlin's list} given in \cite{Mac}, we mention some structures satisfying SAP.
By Theorem \ref{T9050} their posets of copies contain almost disjoint families and maximal antichains of size continuum.\\[-6mm]
\begin{itemize}\itemsep=-1.5mm
\item[-] The posets $\A _\o$, $\B_n$, for $n\le\o$, and $\D$ from Schmerl's list (see Theorem \ref{schmerl});
\item[-] All countable ultrahomogeneous tournaments ({\it Lachlan's list} \cite{Lach}):
$\Q$;
the random tournament, $\T ^\infty$;
the circular tournament, $S(2)$; (see \cite{Che2}, p.\ 18);
\item[-] All {\it Henson's digraphs} with forbidden sets of tournaments \cite{Hen};
(\cite{Mac}, p.\ 11);
\item[-] Digraphs $\Gamma _n$, $n>1$, where $\Gamma _n$ is the Fra\"{\i}ss\'e limit of the amalgamation class of
all finite digraphs not embedding the empty digraph of size $n$;
\item[-] Two ``sporadic" primitive digraphs $S(3)$ and $\CP (3)$;
\item[-] The digraphs $n\ast I_\infty$, for
$2\leq n \leq \o$, which are universal subject to the constraint that non-relatedness is an
equivalence relation with $n$ classes.
\end{itemize}
We remark that some of these structures are not indivisible (so Theorem 5.3 of \cite{Ktow} can not be applied)
for example: $S(2)$, $S(3)$, $\B_n$ and $\C_n$, for $1<n<\o$.
\end{ex}
\section{Ultrahomogeneous partial orders}\label{posets}
Here we prove Theorem \ref{T9092} showing that there are maximal antichains of copies of size $\mathfrak c$
for all ultrahomogeneous partial orders and that $\mathbb A_{\omega}$ and $\mathbb B_{\omega}$ are the only structures on Schmerl's list, for which there are no countable maximal antichains of copies.

First, since the poset $\mathbb P(\mathbb A_{\omega})$ is isomorphic to the poset $\seq{[\o]^{\o},\subset}$, it contains maximal antichains of size $\mathfrak c$,
but does not contain countable maximal antichains.
\subsection{The posets $\B _n$ (disjoint copies of the rational line)}
It is evident that for each $n\le \o$ the poset $\mathbb B_n$ is strongly inexhaustible. So, by Theorem \ref{T9044}, the structure $\mathbb B_n$ satisfies the SAP
and, by Theorem \ref{T9050}, its poset of copies, $\P(\B_n)$, contains maximal antichains of size $\mathfrak c$.
Here we show that, concerning the existence of countable maximal antichains of copies, the finite unions $\mathbb B_n$, $n\in \N$, and
the infinite union $\B_{\o}$ are different. First, the basic case is $\mathbb B_1\cong\Q$.
\begin{lem}\label{T9098}
The poset $\P(\Q)$ contains both $\o$-sized and $\mathfrak c$-sized maximal antichains. 
\end{lem}
\dok
Clearly, the family $\mathcal A=\{I_n:n\in \Z\}$ of the open intervals in $\Q$ given by $I_n=((2n-1)\sqrt{2},(2n+1)\sqrt2)\cap \Q$, for $n\in \Z$, is an antichain in $\P(\Q)$.
If $C\in \P(\Q)$, then $|C\cap I_n|>1$, for some $n\in \Z$, (otherwise we would have $C\hookrightarrow \Z$).
Thus, if $x,y\in C\cap I_n$ and $x<y$, then (since $C\cong \Q$) we have $(x,y)_C\in P(\Q)$ and  $(x,y)_C\subset I_n \cap C$.
So, $\mathcal A$ is a countable maximal antichain in $\P(\Q)$.
\kdok
The following, more general consideration will be used in our analysis of the poset $\P(\B_{\o})$.
For each $i\in \o$, let $\P _i=\la P_i,\leq _i\ra $ be a partial order with a minimum $0_i$ and let $|P_i|\geq 2$.
By $\prod_{i\in \o}\P_i$ we denote the direct product of $\P _i$'s, the poset $\la P,\leq\ra$, where
$P:=\prod_{i\in \o}P_i$ and $\la x_i\ra \leq \la y_i\ra$ iff $ x_i\leq _i y_i$, for all $i\in \o$.
Defining the support of an element $x=\la x_i\ra \in P$ by $\supp(x):=\set{i\in \o:x_i\neq 0_i}$
we consider the suborder $P^{\rm cs}:=\{ x\in P : |\supp (x)|=\o\}$ of the product $\prod_{i\in \o}\P_i$,
call it the {\it countable support product} of $\P_i$'s and denote it by $\prod^{\rm cs}_{i\in \o}\P_i$.
\begin{lem}\label{T9081}
$\prod^{\rm cs}_{i\in \o}\P_i$ does not contain countable maximal antichains.
\end{lem}
\dok
Let $A=\{a^n:n\in \o\}$ be an antichain in $P^{\rm cs}$, where $a^n=\la a^n_i\ra$.
First we show that for different $m,n\in \o$ the set
$$
K_{m,n}=\{i\in \supp(a^m)\cap \supp(a^n):\exists b_i\in P_i \;\; 0_i <_i b_i \leq _i a^m_i, a^n_i\}
$$
is finite. Otherwise, defining $c_i=b_i$, for $i\in K_{m,n}$ and $c_i=0_i$, for $i\in \o \setminus K_{m,n}$, we would have
$c=\la c_i\ra\in P^{\rm cs}$ and $c\leq a^m ,a^n$, which is false.

Let $\la i_n:n\in \o\ra$  be the sequence in $\o$ defined by $i_0=\min(\supp(a^0))$ and
\begin{equation}\label{EQ9001}
i_n=\min \Big(\supp(a^n)\setminus(\{i_0 ,\dots ,i_{n-1}\}\cup K_{0,n}\cup \cdots\cup K_{n-1,n})\Big).
\end{equation}
Let $c=\la c_i:i\in \o\ra$, where $c_{i_n}=a^n_{i_n}$, for $n\in\o$, and $c_i=0_i$, if $i\not\in \{i_n:n\in\o\}$.
For $n\in\o$ we have $i_n\in \supp(a^n)$; thus, $c_{i_n}=a^n_{i_n}>0_{i_n}$ and, hence, $c\in P^{\rm cs}$.

Assuming that $ d\leq a^m,c$, for some $d=\la d_i\ra \in P^{\rm cs}$ and $m\in \o$,
we would have $\supp (d)\subset \supp (c)$ and, hence, $\supp (d)=\{i_n:n\in M\}$, where $M\in [\o]^\o$.
For each $n\in M$ we would have $0_{i_n}<d_{i_n}\leq a^m_{i_n}, a^n_{i_n}$ and, hence, $i_n\in K_{m,n}$,
which is, by (\ref{EQ9001}), impossible for $n>m$.
Thus the element $c$ of $P^{\rm cs}$ is incompatible with all the elements of $A$ and, hence, $A$ is not a maximal antichain in $\prod^{\rm cs}_{i\in \o}\P_i$.
\hfill $\Box$
\begin{te}\label{T9088}
(a) For each $n\in\N$ there is a countable maximal antichain in $\P(\B_n)$.

(b) Each infinite maximal antichain in $\P(\B_{\o})$ is uncountable.
\end{te}
\dok
(a) It is evident that $\mathbb P(\mathbb B_n)=\set{\bigcup_{i<n}\set{i}\times C_i:\forall i<n \;\; C_i\in \mathbb P(\mathbb Q)}$
(see the proof of Theorem 5.1 of \cite{KK1}), which implies that $\mathbb P(\mathbb B_n)\cong \mathbb P(\mathbb Q)^n$.
By Lemma \ref{T9098} there is a countable maximal antichain $A=\{ A_j :j\in \o\}$ in $\P(\mathbb Q)$ and, defining
$\bar A_j :=\la A_j, \Q, \dots ,\Q \ra\in \mathbb P(\mathbb Q)^n$, for $j\in \o$, we obtain a countable antichain $\CA:=\{\bar A_j :j\in \o\}$
in the product $\P(\mathbb Q)^n$. Now, if $\bar C:=\la C_0,\dots ,C_{n-1}\ra\in \mathbb P(\mathbb Q)^n$, then, by the maximality of $A$,
there are $j\in \o$ and $C\in \P(\mathbb Q)$ such that $C\subset C_0 \cap A_j$ and for
$\bar D:=\la C, C_1,\dots ,C_{n-1}\ra$ in the poset $\mathbb P(\mathbb Q)^n$ we have $\bar D \leq \bar C$ and $\bar D \leq \bar A_j$.
Thus $\CA$ is a maximal antichain in the product $\mathbb P(\mathbb Q)^n$.

(b) It is easy to see that the copies of $\mathbb B_\o$ are of the form $\bigcup_{i\in S}\set{i}\times C_i$,
where $S\in [\o]^\o$ and $C_i\in \mathbb P(\mathbb Q)$, for all $i\in S$ (see the proof of Theorem 5.2 of \cite{KK1}).
Thus, the poset $\P(\B_{\o})$ is isomorphic to the countable support product $\prod^{\rm cs}_{i\in \o}\P_i$, where
$\P _i =\la \mathbb P(\mathbb Q)\cup \{ \emptyset \}, \subset\ra$, for all $i\in \o$,
and we apply Lemma \ref{T9081}.
\hfill $\Box$
\subsection{The posets $\C _n$ (dense antichains)}
It is easy to check (see, for example, \cite{KK1}, p.\ 96) that
\begin{eqnarray}
\P (\C _n) & := &\{ n\times A :A\in \P (\Q )\},  \mbox{ for }\;n<\o,  \label{EQ9002}\\
\P (\C _\o) & := & \textstyle\{ \bigcup _{q\in A }C_q \times \{ q \}: A\in \P (\Q ) \land \forall q\in A \;\; C_q \in [\o ]^{\o }\} \label{EQ9003}.
\end{eqnarray}
For $n\leq \o$ and $Z\subset n\times \Q$,  let $\supp(Z)=\{ q\in \Q: Z \cap (n\times \{q\})\neq \emptyset \}$.
Notice that $Z\in \mathbb P(\C_n)$ implies $\supp(Z)\cong \Q$.
\begin{te}\label{T9099}
For each $n$ satisfying $1\leq n\leq \o$ we have

(a) If $\mathcal A$ is a maximal antichain in $\P(\Q)$, then $\mathcal B=\set{n\times A:A \in \mathcal A}$ is a maximal antichain in $\P(\C_n)$
and $|\mathcal A|=|\mathcal B|$;

(b) The poset $\P(\C_n)$ contains both $\o$-sized and $\mathfrak c$-sized maximal antichains.
\end{te}
\dok
(a) Clearly we have $\mathcal B\subset \P(\C_n)$. First we prove that $\mathcal B$ is an antichain.
Assuming that for different $A,A' \in \mathcal A$ there is $Z\in \mathbb P(\mathbb C_n)$
such that $Z\subset n\times A , n\times A'$, we would have $ Z\subset n\times (A\cap A')$ and $\Q \cong\supp (Z)\subset A \cap A'$,
which is impossible since $\mathcal A$ is an antichain in $\P(\Q)$.

Second we prove that $\mathcal B$ is a maximal antichain in $\P(\C_n)$.
If $Z\in \P(\C_n)$, then $\supp(Z)\cong \Q$ and, by the maximality of $\mathcal A$, there are $A\in \mathcal A$ and $B\in \mathbb P(\Q)$
such that $B\subset A\cap \supp(Z)$.
Now, for $Y=\bigcup_{q\in B}Z \cap (n\times \{q\})$ we have $Y\subset Z \cap (n\times A)$
and $Y\in \P (\C_n)$ because for each $q\in \supp(Z)$ there is a bijection between $n$ and $Z \cap (n\times \{q\})$.
So, $Y$ witnesses the compatibility of $Z$ and $n\times A\in \mathcal B$.

(b) follows from (a) and Lemma \ref{T9098}.
\hfill $\Box$
\subsection{The poset $\D$ (the random poset)}
Recall that $\mathbb D=\la D,<\ra$ is the unique, up to isomorphism,
countable ultrahomogeneous partial order which embeds all countable partial orders.
Since the structure $\mathbb D$ satisfies the SAP,
by Theorem \ref{T9050} the poset $\P(\mathbb D)$ contains maximal antichains of size $\mathfrak c$.
So, Theorem \ref{T9091} given below completes the proof of Theorem \ref{T9092}.

First we recall some definitions and facts from \cite{KK1}
which will be used in the proof of Theorem \ref{T9091}
(see  Fact 3.1, Fact 3.2 and Lemma 3.3 of \cite{KK1})
and note that $\parallel$ will denote the {\it incomparability relation}:
$p\parallel q \Leftrightarrow p\neq q \land \neg p< q \land \neg q< p$.

\begin{df}\rm\label{D5000}
Let $\mathbb P= \langle P, < \rangle$ be a partial order. By
$C(\mathbb P )$ we denote the set of all triples $\langle L,G,U
\rangle$ of pairwise disjoint finite subsets of $P$ such that:

(C1) $\forall l\in L \;\; \forall g\in G \;\; l<g$,

(C2) $\forall u\in U \; \forall l\in L \; \neg u<l $,

(C3) $\forall u\in U \; \forall g\in G \; \neg g<u $.

\noindent For $\langle L,G,U \rangle \in C(\mathbb P )$, let $P
_{\langle L,G,U \rangle}$ be the set of all $p\in P\setminus (L
\cup G\cup U)$ satisfying:

(S1) $\forall l\in L \;\; p>l $,

(S2) $\forall g\in G \; p<g $,

(S3) $\forall u\in U \; p \parallel u $.
\end{df}
\begin{fac}\label{T5040}
A countable partial order $\mathbb P = \langle P, < \rangle$ is (isomorphic to) a countable random poset iff
$P_{\langle L, G ,U \rangle}\neq \emptyset$, for each $\langle L, G ,U \rangle \in C(\mathbb P )$.
\end{fac}
\begin{fac}\label{T5003}
Let $\mathbb P= \langle P, < \rangle$ be a partial order and
$\emptyset \neq A\subset P$. Then

(a) $C(A, < )= \{ \langle L,G,U \rangle \in C(\mathbb P ): L,G,U
\subset A\}$;

(b) $A_{\langle L,G,U \rangle}= P_{\langle L,G,U \rangle}\cap A$,
for each $\langle L,G,U \rangle \in C(A, < )$.
\end{fac}
\begin{fac}\label{T5006}
Let $\mathbb D = \langle D, < \rangle$ be a countable random
poset. Then

(a) $D_{\langle L, G ,U \rangle}\in \mathbb P (\mathbb D )$, for each $\langle L, G ,U \rangle \in C(\mathbb D )$;

(b) If $C\subset D$ and $A\not\subset C$, for each $A\in \mathbb P (\mathbb D)$, then $D\setminus C\in \mathbb P (\mathbb D )$;

(c) If $\mathcal L\subset \mathbb P (\mathbb D)$ is a chain, then
$\bigcup\mathcal L\in \mathbb P (\mathbb D)$.
\end{fac}
\begin{te}\label{T9091}
There is a countable maximal antichain in $\P(\mathbb D)$.
\end{te}
\dok
Let ${\mathcal C}$ be a maximal chain in $\D$. Assuming that $x=\max {\mathcal C}$ (resp.\ $x=\min {\mathcal C}$)
we would have $\D_{\la \{ x\},\emptyset ,\emptyset \ra}=\emptyset $ (resp.\ $\D_{\la \emptyset , \{ x\},\emptyset \ra}=\emptyset $).
Thus ${\mathcal C}$ is an unbounded chain in $\D$ and, since $|{\mathcal C}|=\o$,
it has a cofinal subset isomorphic to $\o$ and a coinitial subset isomorphic to $\o ^*$.
This implies that $\D$ contains an unbounded copy of the integers.
W.l.o.g.\ we suppose that $\Z$ itself is that copy. For $m\in \mathbb Z$, let
$A_m=\{x\in \mathbb D: x< m \}$ and $X_m=A_m\setminus (A_{m-1}\cup \{m-1\})$; that is,
\begin{equation}\label{EQ9004}
X_m=\{ x\in D : x<m \land x\not\leq m-1\}.
\end{equation}
{\it Claim.}
$X_{m}\in \P(\mathbb D)$, for every $m\in \Z$.\\[-2mm]

\noindent
{\it Proof of Claim.}
We show that $X_m\cong \D$ using Fact \ref{T5040}.
If $\la L,G,U \ra \in C(X_m)$, then by Fact \ref{T5003}(a) we have $\la L,G,U \ra \in C(\D)$ and $L,G,U \subset X_m$.
By Fact \ref{T5003}(b), we have to show that $D_{\seq{L,G,U}}\cap X_m\neq \emptyset$. There are four cases.

{\it Case I}: $L\neq\ps$ and $G\neq\ps$.
Since $\la L,G,U \ra \in C(\D)$, there is $d\in \mathbb D_{\seq{L,G,U}}$
and we show that $d\in X_m$.
Since $\emptyset \neq G \subset X_m$, for $g\in G$  we have $d<g<m$.
Assuming that $d\leq m-1$,
for $l\in L$ we would have $l<d \leq m-1$, which is false because $L\subset X_m$ and, by (\ref{EQ9004}), $l\not\leq m-1$.
So, $d\not\leq m-1$ and, by (\ref{EQ9004}), $d\in X_m$.

{\it Case II}: $L=\ps$ and $G\neq \ps$.
Suppose that $\la \ps,G,U\cup \{m-1\}\ra\not\in C(\D)$.
Then (C3) fails
and, since $\la \ps,G,U \ra \in C(\D)$, there is $g\in G$ such that $g<m-1$.
But, since $G\subset X_m$, this is impossible by (\ref{EQ9004}).
Thus $\la\ps,G,U\cup\{m-1\}\ra\in C(\D)$
and, hence, there is $d\in \D_{\la \ps,G,U\cup\{m-1\}\ra}$, which implies $d\in \D_{\la \ps,G,U\ra}$.
We prove that $d\in X_m$.
First, for $g\in G$ we have $d<g<m$.
Second, since $d\parallel m-1$, we have $d\not\leq m-1$
and, by (\ref{EQ9004}), $d\in X_m$ indeed.

{\it Case III}: $L\neq\ps$ and $G=\ps$.
Consider the condition $\seq{L,\set{m},U}$.
Since $L\subset X_m$, by (\ref{EQ9004})  we have $ l<m$, for all $l\in L$, and (C1) is true.
(C2) is true because $\la L,\emptyset ,U \ra \in C(\D)$.
Since $U\subset X_m$, by (\ref{EQ9004}) for each $u\in U$ we have $\neg m<u$ and (C3) is true.
So, $\seq{L,\set{m},U}\in C(\D)$,
there is $d\in \D_{\seq{L,\set{m},U}}\subset \D_{\seq{L,\ps,U}}$ and we show that $d\in X_m$.
Clearly $d<m$ and $d\leq m-1$ would imply that $l< d\leq m-1$, for some $l\in L$,
which is false because $\emptyset \neq L\subset X_m$.
So, $d\not\leq m-1$ and $d\in X_m$.

{\it Case IV}: $L=\ps$ and $G=\ps$.
Suppose that $\seq{\ps,\set{m},U\cup\set{m-1}}\not\in C(\D)$.
Then (C3) fails so there is $u\in U\cup\set{m-1}$ such that $m<u$.
Since $m\not<m-1$ we have $u\in U$, which is false because $u\in X_m$.
Thus $\seq{\ps,\set{m},U\cup\set{m-1}}\in C(\D)$,
there is $d\in \D_{\seq{\ps ,\set{m},U\cup\set{m-1}}}\subset \D_{\seq{\ps,\ps,U}}$
and we show that $d\in X_m$.
Clearly $d<m$
and $d\parallel m-1$ gives $d\not\leq m-1$;
so,  $d\in X_m$ indeed.
\kdok
\noindent
{\it Claim.}
$X=\bigcup_{m\in \mathbb Z}X_m$ is isomorphic to the random poset.\\[-2mm]

\noindent
{\it Proof of Claim.}
By Fact \ref{T5006}(a) we have $A_m=\D_{\seq{\ps,\{m\},\ps}}\in \P(\D)$, for all $m\in \Z$.
Since $\{ A_m :m\in \Z\}$ is a chain in $\P(\D)$, Fact \ref{T5006}(c) gives $A=\bigcup_{m\in \Z}A_m\in \mathbb P(\mathbb D)$
Since for each $m\in \mathbb Z$ we have $X_m\subset A_m$ and $X_m\cap \Z=\emptyset $ it follows that $X\subset A\setminus \Z$.
Conversely, since $\bigcap _{m\in \Z}A_m=\emptyset$ (by the unboundedness of $\Z$),
for $x\in A\setminus \Z$ there is $m\in \Z$ such that $x\in A_m \setminus A_{m-1}$
and, since $x\neq m-1$, we have $x\in X_m \subset X$; so $X= A\setminus \Z$.
Since $A\cong \D$ and $\Z$ does not contain copies of $\D$, by Fact \ref{T5006}(b) we have $A\setminus \Z\in \mathbb P (\mathbb D )$
that is, $X\in\mathbb P(\D)$.
\kdok
Finally, we prove that $\mathcal A=\set{X_m:m\in \Z}$ is a maximal antichain in $\P(X)$.
Since $X_m \cap X_n=\emptyset$, for different $m,n\in \Z$,  $\mathcal A$ is an antichain in $\P(X)$. For a proof of its maximality
we take $C\in \P(X)$ and first, towards a contradiction, suppose that
\begin{equation}\label{EQ9005}
\forall m\in \Z \;\;\forall x,y \in C \cap X_m\;\;(x\neq y \Rightarrow x\parallel y).
\end{equation}
Let us fix $x\in C$ and $m_0\in \Z$, where $x\in C\cap X_{m_0}$.
Since $C_{\la \{ x \}, \emptyset , \emptyset\ra}\neq \emptyset$ there are $m\geq m_0$ and $z\in C\cap X_m$ such that $x<z$,
which by (\ref{EQ9005}) implies that $m>m_0$.
Thus the set $M:=\set{m>m_0: \exists z\in C\cap X_m\ \ x<z}$ is non-empty.
Let $m_1=\min M$ and let us pick $z\in C\cap X_{m_1}$ such that $x<z$.
Then, since  $C_{\la \{ x \}, \{ z \} , \emptyset\ra}\neq \emptyset$
there are $m\in \Z$ and $y\in C\cap X_m$ such that $x<y<z$.
Now, $m\in M$, we have $m_0\leq m \leq m_1$,
and, by (\ref{EQ9005}),  $m_0< m < m_1$,
which is impossible because  $m_1=\min M$.

Thus (\ref{EQ9005}) is false and, hence,
there are $m\in \Z$ and $x,y\in X_{m}\cap C$ such that $x<y$.
Since $C\in \mathbb P(X)$, by Fact \ref{T5006}(a) we have $C_{\seq{\set{x},\set{y},\ps}}\in \mathbb P(X)$.
If $t\in C_{\seq{\set{x},\set{y},\ps}}$, then $x<t<y<m$
and $t\leq m-1$ would imply $x<m-1$, which is false because $x\in X_m$;
so, $t\in X_m$. Thus, $C_{\seq{\set{x},\set{y},\ps}}\subset C\cap X_m$,
that is, $C$ and $X_m$ are compatible elements of $\P (X)$
and $\mathcal A$ is a maximal antichain in $\P (X)$.
\hfill $\Box$

\end{document}